\documentclass[journal]{IEEEtran}
\usepackage{ifpdf}

\usepackage{algorithm}
\usepackage{algpseudocode}
\usepackage{color}
\usepackage{url}
\ifpdf
\usepackage[pdftex]{graphicx}
\else
\usepackage[dvips]{graphicx}
\fi
\usepackage{tabularx}
\DeclareGraphicsExtensions{.eps}
\usepackage[caption=false,font=footnotesize]{subfig}
\usepackage{amsmath,amsthm,amsfonts,amssymb}
\usepackage[noadjust]{cite}
\usepackage{nomencl}
\usepackage{enumerate}
\usepackage{bm}
\theoremstyle{definition}
\newcommand{\ud}[1]{_\mathrm{#1}}
\newcommand{\up}[1]{^\mathrm{#1}}

\newtheorem{theorem}{\textbf{Theorem}}

\newtheorem{remark}[theorem]{\textbf{Remark}}

\newcommand{\rev}[1]{\textcolor{black}{#1}}
\newcommand{\revv}[1]{\textcolor{black}{#1}}

\begin{document}


\title{Dynamic Valuation of Battery Lifetime}

\author{Bolun~Xu,~\IEEEmembership{Member,~IEEE}
	\thanks{B.~Xu is with Columbia University, USA (emails: bx2177@columbia.edu). }
}  

\maketitle

\begin{abstract}
This paper proposes a dynamic valuation framework to determine the opportunity value of battery capacity degradation in grid applications based on the internal degradation mechanism and utilization scenarios. The proposed framework follows a dynamic programming approach and includes a piecewise linear value function approximation solution  that  solves the optimization problem over a long planning horizon. The paper provides two case studies on price arbitrage and frequency regulation using real market and system data to demonstrate the broad applicability of the proposed framework.  Results show that the battery lifetime value is critically dependent on both the external market environment and its internal state of health. \rev{On the grid service side, results show that second-life batteries can provide more than 50\% of the value compared to new batteries, and frequency regulation provides two times more revenue than price arbitrage throughout the battery lifetime.}

\end{abstract}




\section{Introduction}

The global lithium-ion battery production capacity is on track to surpass 1000~GWh per year before 2023, while the pack price has decreased by more than 80\% over the past ten years and is now approaching \$100/kWh~\cite{blochbreakthrough}. However, rapid technology advancements and cost reductions may not be sufficient for economic and sustainable battery deployments in grid applications. One must characterize the cost of operating the storage, including efficiency, degradation, and future opportunities. In power systems, operators will schedule dispatch based on  bids from various resources, including traditional thermal generators, renewables including wind and solar, and batteries.  In an efficient and adequately regulated market, bids will reflect the true operating cost of these resources. The market-clearing result will optimally be allocated resources that maximize social welfare~\cite{david2001market}. In the example of conventional thermal generators, this operating cost should be primarily based on the fuel cost and the heat rate efficiency of the generator.  

Yet, quantifying the capacity value and operating cost for battery energy storage is non-trivial. \rev{Batteries have negligible immediate cycling costs as they consume no fuel to store or generate energy.} However, cycling a battery leads to material and structural degradation inside the battery cell~\cite{vetter2005ageing}, permanently reducing the storage capacity and battery lifetime. For example, recent lab tests show one full cycle will reduce an NMC (lithium-manganese-cobalt-oxide) cell capacity by about 0.04\%, meaning this NMC cell may only be cycled 500 times before reaching its \rev{rated end-of-life}~\cite{preger2020degradation}. The capacity reduction accumulates throughout the entire remaining battery lifetime. \rev{As a comparison, the start-up and shut-down costs of thermal generators have been incorporated into unit commitments as they increase the maintenance cost significantly~\cite{kumar2012power}. In contrast, the new pay-for-performance regulation market has incorporated a mileage payment to compensate generators for their frequent generation set-point cycles to follow regulation signals~\cite{xu2016comparison}.}
To this end, battery cycles at the cost of its lifetime and reduce future opportunities. Hence the cost of cycling must be systematically incorporated into economic plannings.

This paper proposes a novel dynamic programming approach that accurately factors the opportunity value of battery degradation accounting cycle degradation, calendar degradation, discount rate, and the application utility. The contribution in this paper is listed as follows:
\begin{enumerate}
    \item It proposes a novel and generalized dynamic programming formulation to model the impact of battery degradation over multi-year operation planning.
    \item It develops two solution algorithms for solving the dynamic valuation framework aiming at different application scenarios: a piecewise linear approximation approach for optimization and a simulation-based algorithm for stochastic control and simulations.
    \item It demonstrates the proposed valuation framework in energy price arbitrage and frequency regulation using real market data over a multi-year project lifetime.
    \item It compares the performance of new and second-life batteries and concludes that second-life batteries are an economical choice for providing grid services.
\end{enumerate}

The rest of the paper is organized as follows. Section II includes literature reviews. Section III proposes the formulation and Section IV introduces the solution methods. Section V introduces the case study settings and Section VI presents the results. Section VII concludes the paper.

\section{Background and Literature Review}


Battery cells are constantly losing their usable energy capacity due to degradation,  eventually causing them to retire and be recycled~\cite{vetter2005ageing, zakeri2015electrical}. In consumer applications such as personal electronics or electric vehicles, users do not change their use pattern based on  battery degradation, and battery lifetime must meet the expected device application lifetime. For example, most EV manufacturers are offerings 8-year/100,000 mile warranties for their batteries~\cite{doe_ev} which ensures batteries in the EV will function properly and meet a specified mileage range within the period and distance in the warranty.

However, the consideration of battery degradation can be much more sophisticated in grid-interactive applications since battery owners and operators have the decision flexibility on whether or not to use,  how much to use, or which battery to use among a fleet of batteries at any period. From the perspective of battery participants in electricity markets who wish to maximize the investment profit, the use of batteries should be optimized considering the limited lifetime and the volatile market conditions. While from the perspective of power system dispatch, battery degradation should be incorporated into economic dispatch to minimize the social cost and environmental impacts, including battery production, recycling, and disposal. In both perspectives, the operator must trade-off short-term battery utilization with the long-term life expectancy. 


Modeling battery degradation mechanisms in power system optimizations have received much attention over the past few years, and many degradation models have been incorporated into optimization in computation tractable formulations. These include the most basic energy throughput model that assumes capacity fade is proportional to the accumulated energy charged and discharged from the battery~\cite{xu2017scalable}; the C-rate model, which considers increased degradation rate from high charge and discharge power~\cite{wang2019power}; and the Rainflow-based cycle depth model which considers the nonlinear degradation impact from cycles of various depths~\cite{xu2017factoring, shi2017using, he2015optimal}. More sophisticated degradation models, including  equivalent circuit models and  single-particle models, have also been explored in optimization, but the computation was not tractable~\cite{reniers2018improving}. Battery degradation models have been incorporated into optimization constraints. Still, the operator must first specify the maximum allowable degradation  during the operation horizon~\cite{mohsenian2015optimal}, which conserves the battery utilization under power system operation uncertainties. Hence, it is more desired to incorporate degradation models into the objective function as part of the battery operating cost.

Yet, there is one more step before including a degradation model in the optimization problem as a cost term: an incremental cost of degradation must be specified. For example, what should be the cost if the considered battery lost 0.01\% of usable energy capacity? A common approach in prior studies and industry practices is to amortize the battery replacement pack cost and prorated it to marginal capacity losses, under the intuition that if the current battery pack reached end-of-life, a new pack could be purchased and installed to continue the operation of the battery energy storage system. This replacement cost is used in instead of the operating cost. For example, the amortized cost of a full cycle could be around \$100/MWh if the battery pack cost is \$200/kWh, assuming a 2000 rated cycle life. Cost amortization is a convenient approach, as the purchase cost or replacement cost is usually known to the battery operator. It does not add computation complexities but also comes with a drawback: it does not maximize the profit or utilization of the battery throughout the lifetime and may not even ensure the revenue pay-off the battery pack cost. A battery price arbitrage case study~\cite{xu2017factoring} shows the cost amortization approach can improve the lifetime battery market revenue by 200\% compared to ignoring the degradation model, but still, the resulting lifetime NPV (net present value) revenue for the battery is about \$200/kWh which is less than the \$300/kWh pack cost used in the cost amortization.

The drawback of cost amortization was first tackled by He et al.~\cite{he2018intertemporal,he2020power}, who stated the battery degradation cost should be factored based on the loss of future opportunity instead of the pack cost and developed a computation framework for estimating the marginal cost of battery degradation that maximizes the lifetime battery utilization. Still, the marginal cost was modeled as a static value, and its dependency over the battery state-of-life was not considered. Degradation permanently reduces the usable battery capacity, losses occurred at a newer battery have different consequences than at a more degraded battery. On the other hand,  prior degradation models  neglect calendar degradation in operations and planning, but doing so will significantly over estimate the cost of battery degradation. For example, a high-cost battery pack will not cycle much because its amortized degradation cost is much higher than market payments, thus this battery will ``waste'' most of its capacity due to calendar degradation. Therefore, a systematic battery degradation valuation framework must be designed to incorporate calendar degradation and dependency over state-of-life.

\section{Formulation}

\begin{figure}[t]
    \centering
    \includegraphics[trim = 45mm 70mm 90mm 35mm, clip, width = .95\columnwidth]{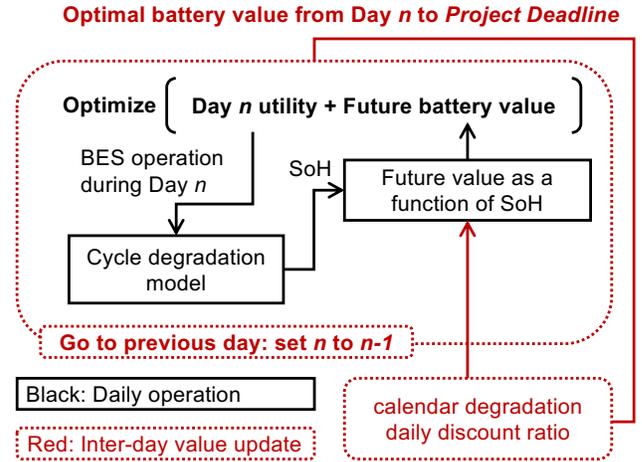}
    \caption{The proposed dynamic battery valuation framework. This framework works backward and keeps updating a value function with respect to battery SoH. This function represents the value of the remaining battery SoH from the end of the current operating day till the project deadline, at which the battery has no more arbitrage value. At each operating day, the battery is optimized according to daily revenue income and the change in the battery SoH, based on a cycle degradation model. The result of this optimization is assembled into the new battery future value function, and this framework works recursively backward to evaluate the battery lifetime across the entire project duration.}
    \label{fig:framework}
\end{figure}

Fig.~\ref{fig:framework} illustrates the proposed valuation framework. The key idea is to define the value of the battery as a time-variant function of the battery state-of-health (SoH). \rev{The calculated value function represents the maximized value of the battery at a given SoH from the current day until the project deadline, spanning from a new battery with 100\% SoH to its end-of-life (EoL):
\begin{itemize}
    \item \textbf{State-of-Health (SoH)} is the ratio between the current remaining battery capacity and the rated new battery capacity. For example, a new battery has 100\% SoH, while a battery with 80\% SoH means it has lost 20\% of the energy capacity due to degradation.
    \item \textbf{End-of-Life (EoL)} is the SoH at which the battery can no longer serve its application purposes and must be retired from the battery system, often due to increased impedance and reduced thermal stability~\cite{vetter2005ageing}. Many  manufacturers rate their lithium-ion batteries with an 80\% EoL, but many batteries can function afterward, reaching around 60\% SoH~\cite{birkl2017degradation}.
\end{itemize}
}


The valuation framework follows a standard dynamic programming approach~\cite{bellman2015applied}. It initializes the battery value at the end of the project period, either as zero or a resale value, assuming the battery project terminates after the project lifetime. Any battery cells left will have no more value at the project site. The framework then works backward in time and recursively updates the opportunity value of battery degradation based on the daily revenue. \rev{For example, assume the battery owner can cycle the battery at day $n$ to make \$200 but reduce the battery lifetime from 90\% to 89.5\% SoH. Since the algorithm works backward, we know from the value function that a 90\% SoH battery can make a total NPV (net present value)  revenue of \$10,000  from day $n+1$ until the end of the battery or the project lifetime. In comparison, an 89.5\% SoH battery can only make  an NPV of \$9,900. Hence, reducing SoH from 90\% to 89.5\% has an opportunity cost of \$100, which is lower than the \$200 profit opportunity, indicating the battery owner should cycle the battery. This action will also update the opportunity value of the battery capacity at day $n$ to \$10,200 at 90.5\% SoH and \$10,100 at 90\% SoH.\footnote{For simplicity, this example ignores calendar degradation, daily discount ratio, and the impact of capacity reduction in revenue collection. These factors are addressed in the full formulation.}}

We start with a deterministic formulation follow a dynamic programming approach which is recursively defined as:
\begin{subequations}\label{eq:dbv}
\begin{align}\label{eq:dbv_obj}
    V_{n}(E_{n}) := \max_{\bm{p}_n\in\mathcal{P}(E_n)} O_n(\bm{p}_n) + \gamma V_{n+1}(E)
\end{align}
subjects to the capacity degradation model
\begin{align}
    E = E_{n} - \big[D\ud{cyc}(\bm{p}_n) + D\ud{cal}\big]E^0\label{eq:dbv_c3}
\end{align}
and the battery operation subjects to the following constraints
\begin{align}
    & e_{n,t} = e_{n,t-1} - H(p_{n,t}) \label{eq:dbv_c1}\\
    & 0\leq e_{n,t} \leq E_{n} \label{eq:dbv_c2}\\
    & p_{n,t} \in [-P,P],\quad e_{n,T} \geq e_{n,0} \label{eq:dbv_c4}
\end{align}
\end{subequations}
where
\rev{
\begin{itemize}
    \item $n\in\{1,\dotsc, N\}$ is the index of the valuation stages and are designed as operating days (24 hours) in this problem.
    \item $t\in\{1\dotsc,T\}$ is the index of the operation intervals at which granularity the battery control decisions are updated. 
    \item $\bm{p}_n = \{p_{n,t}\in\mathbb{R}\}$ is the set of the battery dispatch set-points during operating day $n$, positive for discharge and negative for charge.
    \item $E_n$ is the battery energy capacity on day $n$.
    \item $\mathcal{P}(E_n)$ is the feasibility set of the battery operation as described in \eqref{eq:dbv_c1}--\eqref{eq:dbv_c4}. It is dependent on the battery energy capacity $E_n$ to show the operation depends on the remaining battery energy capacity. 
    \item $e_{n,t}$ is the battery stored energy during interval $t$ on day $n$.
    \item $O_n(\cdot)$ is the daily revenue function as a result of the battery operation profile $\bm{p}_n$.
    \item $V_{n+1}(\cdot)$ is the value-to-go function of battery's remaining capacity $E$ describing the remaining value of the battery capacity at the end of day $n$ (the same to at the beginning of day $n+1$) throughout the rest of the project lifetime.
    \item $\gamma$ is the daily discount rate.
    \item $D\ud{cyc}(\bm{p}_n)$ is the cycle degradation model that calculates the reduction of battery capacity as a result of the battery operating profile $\bm{p}_n$.
    \item $D\ud{cal}$ is the daily calendar degradation rate which we assume to be a constant.
    \item $E^0$ is the rated energy capacity of the battery, i.e., without any degradation. 
    \item $H(\cdot)$ models the efficiency of the battery from charging or discharging. For example, if we assume a linear charging and discharging efficiency rate $\eta$ and battery control is updated at five-minute resolution, $H(\cdot)$ can be written as $([p_{n,t}]^-\eta + [p_{n,t}]^+/\eta)/12$ where $[x]^-$ gets the negative component of $x$ and $[x]^+$ get the positive component of $x$. 
\end{itemize}
}
\eqref{eq:dbv_c3} models the daily change of the battery remaining capacity due to cycle and calendar degradation, in which $D\ud{cyc}(\bm{p}_n)$ is the cycle degradation model that calculates the reduction of battery capacity as a result of the battery operating profile $\bm{p}_n$, and $D\ud{cal}$ is the daily calendar degradation rate which we assume to be a constant. The degradation models are normalized and are multiplied with the rated energy capacity $E^0$ to reflect the change of capacity in MWh. \eqref{eq:dbv_c1} models the battery state-of-charge evolution constraint using the efficiency function $H$. \eqref{eq:dbv_c2} enforces the battery energy capacity limit $E_n$, and \eqref{eq:dbv_c4} models the battery power rating and enforces the battery final state-of-charge at the end of the day must be no smaller than at the beginning of the day.

A key \emph{assumption} we made in \eqref{eq:dbv} is that the battery energy capacity remains constant during the day and the incremental degradation caused by the daily operation is added to the battery at the end of the day. Although this assumption is different from the underlying degradation process, the error is negligible as lithium-ion batteries usually has a lifetime over five to ten years to reach an end-of-life around 70\%, hence the daily capacity error from this assumption is about 0.01\% to 0.03\% of the battery's rated energy capacity. This error is sensitive to the duration of valuation stages, for example, if we change the duration of $n$ from one day to one month, the error becomes 0.3\% to 0.9\%. 

This formulation applies to most battery application scenarios in which battery operation can be formulated as a convex optimization problem. For example, most battery grid applications can be modeled with convex objectives including economic dispatch~\cite{kirschen2018fundamentals}, price arbitrage~\cite{krishnamurthy2017energy}, frequency regulation~\cite{xu2018optimal}, or peak shaving~\cite{shi2017using}. Convex degradation models include the most common energy throughput model~\cite{bishop2013evaluating}, while recent studies have shown more complex Rainflow-based degradation models are also convex~\cite{shi2018convex,xu2017factoring}.

\rev{Notably, non-convex models such as binary constraints can be incorporated into the valuation framework, and still likely to provide tractable computation time since the framework is based on solving daily optimizations or simulations. However, the overall optimality of the valuation result cannot be guaranteed. }
\rev{Degradation models $D\ud{cyc}(\cdot)$ and $D\ud{cal}(\cdot)$ can optimally include other factors if formulate into a convex form. For example, the battery state-of-charge $e_{n,t}$ may also be incorporated into the cycle and calendar degradation model according to lab test results~\cite{ecker2014calendar}. The battery cell temperature can also be included as a parameter in the degradation functions. However for simplicity of the mathematical presentation in this paper, we will not include an exhaustive list of degradation factors. }


\section{Solution Method}
We develop two algorithms for solving the dynamic valuation problem, both algorithms use a piecewise linear value function approximation approach. The first one, Algorithm 1, is a deterministic optimization model that integrates value function approximation with the decision making process. The second one, Algorithm 2, is a simulation-based version which is more suitable for stochastic control scenarios or simulations.

\subsection{piecewise Linear Value Function Approximation}
We propose a solution method for the dynamic valuation problem \eqref{eq:dbv} by constructing a piecewise linear concave approximation of the value function $V_n$. Note that if the utilization function $O_n$ is concave and the cycle degradation function $D\ud{cyc}$ is convex, then \eqref{eq:dbv} is a concave optimization problem and the resulting value function $V_n$ is therefore a concave function that can be approximated using a piecewise linear approach. 

We initialize a set of $i\in\{1,\dotsc,I\}$ monotonically decreasing battery energy capacity sampling points $\{E\up{s}_1,\dotsc, E\up{s}_I\}$, with $E\up{s}_1$ be the rated energy capacity of the battery with no degradation, and $E\up{s}_I$ be the end-of-life battery capacity, so the sample covers the entire battery capacity range. \revv{A higher number of segments provides better accuracy but will also increase the computation complexity. Because the degradation rate is relatively smooth in most grid applications~\cite{xu2016modeling}, a 1\% discretization granularity/sensitivity should provide reasonably good approximations.}

During operating day $n$, we repetitively solve \eqref{eq:dbv} using each initial battery capacity sample and obtain the resulting optimal objective value $v\up{s}_{i,n}$, with the only exception being $E\up{s}_I$ whose corresponding value should be kept to zero, i.e. $v\up{s}_{I,n} = 0$ for all $i\in\{1,\dotsc,I\}$, indicating a battery that reached end-of-life has zero future value. Therefore $V_n$ can be approximated using $(E\up{s}_i, v\up{s}_{i,n})$ as the following concave piecewise linear approximation
\begin{align}\label{eq:pwl}
    V_n(E) \approx \max_{v} \Big\{& v | v \leq v\up{s}_{i,n} + (v\up{s}_{i,n} - v\up{s}_{i+1, n}) \frac{E - E\up{s}_i}{E\up{s}_i - E\up{s}_{i+1}},\nonumber\\
    &\forall i\in\{1,\dotsc, I-1\} \Big\}
\end{align}
which can be directly added to \eqref{eq:dbv} to replace $V_n(E)$ and the maximization can be unified with the maximization operation in \eqref{eq:dbv} in which the objective function becomes
\begin{subequations}\label{eq:p2}
\begin{align}\label{eq:p2_obj}
    v\up{s}_{i,n} := \max_{\bm{p}_n} O_n(\bm{p}_n) + \gamma v
\end{align}
subjects to the piecewise linear value function approximation \eqref{eq:pwl} and the degradation constraint which is rewritten as
\begin{align}
    E = E\up{s}_i -  \big[D\ud{cyc}(\bm{p}_n) + D\ud{cal}\big] E^0
\end{align}
and the battery operating constraints \eqref{eq:dbv_c1}--\eqref{eq:dbv_c4}.
\end{subequations}

The full solution algorithm \textbf{Algorithm 1} is listed as 
\begin{enumerate}
    \item Set $N\to n$, start from the last battery project day.
    \item Initialize $0\to v\up{s}_{i,N+1}$ for all $i\in\{1,\dotsc,I\}$, indicating the battery has no future value.
    \item For $i=1$ to $I-1$, solve \eqref{eq:p2} and record the optimized objective value to $v\up{s}_{i,n}$; note $v\up{s}_{I,n}$ is always zero (the value of end-of-life battery).
    \item If $n=1$, return; otherwise, set $n-1\to n$ and go to Step 3 (go to the previous day).
\end{enumerate}
This algorithm has a linear complexity with respect to the number of days or stages considered in the valuation problem, hence the dynamic valuation framework is computation tractable over any planning horizon.



\subsection{\rev{Non-anticipatory Formulation and Algorithm}}

\rev{We extend the dynamic valuation framework to non-anticipatory control scenarios to address consideration of uncertainties or complex battery models. In this framework, we assume the battery power $\bm{p}$ being regularly updated in time sequence in accordance with a stochastic model or a control rule, instead of using a deterministic optimization as described in \eqref{eq:dbv}. We generalize the battery operation as product of a decision policy $G$, the dynamic valuation framework is thus represented below using the piecewise linear value function discretization method as presented in \eqref{eq:pwl}.} This algorithm, \textbf{Algorithm 2}, is listed below
\begin{subequations}
\begin{enumerate}
    \item Calculate the marginal cost of battery degradation $C_{i,n}$ at segment $i$ during operating day $n$ using the value function results from the next operating day $n+1$ as the slope between the neighbouring segments as
    \begin{align}
        C_{i,n} &\leftarrow \frac{v\up{s}_{i,n+1} - v\up{s}_{i+1, n+1}}{E\up{s}_i - E\up{s}_{i+1}}
    \end{align}
    note that it is important that the slope being calculated using the more degraded segment $i+1$ instead of $i-1$ as battery degradation is a unidirectional process: a battery will only loss capacity but not gain more. 
    \item Calculate the battery operation profile $\bm{p}$ using a decision rule $G$ in which the inputs include the uncertainty model $\bm{\varepsilon}$, the marginal degradation cost $C_{i,n}$, and the corresponding energy capacity $E\up{s}_{i}$ 
    \begin{align}
        \bm{p}_{i,n} &\leftarrow G(\bm{\varepsilon}_n, C_{i,n}, E\up{s}_{i}) 
    \end{align}
    here we list $C_{i,n}$ and $E\up{s}_{i}$ as inputs to the decision rule to show $\bm{p}_{i,n}$ is dependent on the battery energy capacity and the marginal degradation cost. Other constraints such as efficiency and power should be included in $G$ but are not listed explicitly as inputs here because we assume they do not vary with the battery lifetime. 
    \item Calculate the resulting degradation from the operation profile $\bm{p}_{i,n}$  as
    \begin{align}
        D_{i,n} &\leftarrow D\ud{cyc}(\bm{p}_{i,n}) + D_{cal} 
    \end{align}
    \item Store the value function result for operating day $n$ as
    \begin{align}
        v\up{s}_{i,n} &\leftarrow \gamma (v\up{s}_{i,n+1}  - C_{i,n}D_{i,n}E^0) +  O_n(\bm{p}_{i,n})
    \end{align}
    $C_{i,n}D_{i,n}E^0$ calculates the change in the battery value as the product of the marginal cost of degradation $C_{i,n}$, the degradation rate $D_{i,n}$, and the rated battery capacity $E^0$. \revv{$v\up{s}_{i,n+1}  - C_{i,n}D_{i,n}E^0$ represents the opportunity value of the remaining battery lifetime in the next day, and is multiplied with the daily discount ratio $\gamma$ in the current day operation as in \eqref{eq:dbv_obj}. }
\end{enumerate}
The battery dynamic valuation is performed by iterating above steps $n$ from $N$ to 1 with an inner iteration of $i$ from $I-1$ to 1 embedded, same as the iteration setting in Algorithm 1.
\end{subequations}

Algorithm 2 is an generalized version of Algorithm 1, i.e., Algorithm 2 provides the same result to Algorithm 1 if replace $G$ with the optimization model as described in \eqref{eq:dbv}. Algorithm 2 decouples the decision making process $G$ from the battery value calculation, so that it accommodates more realistic control settings such as non-anticipatory controls or even simulations. 

\begin{remark}\textbf{Non-convex and nonlinear battery models.}
\rev{Algorithm 2 is more efficient in incorporating highly nonlinear battery simulation and degradation models since the value function update is decoupled from the battery operation, so the battery operation can be performed using simulations instead of having to solve an optimization problem.} 

\begin{remark}\textbf{Piecewise linear granularity.}
The piecewise linear approximation granularity in Algorithm 2 cannot be too fine, since it assumes the marginal cost of battery degradation is a constant throughout the operating day $n$. The piecewise linear segment size in the degradation value function must be larger than the maximum possible incremental degradation from a valuation stage.
\end{remark}

\end{remark}

\section{Case Study}\label{sec:sv}




\subsection{Battery and Degradation Model}
In all case studies we consider a battery energy storage made up with $\mathrm{Li(NiMnCo)O_2}$-based 18650 lithium-ion battery cells with the following parameters unless otherwise specified:
\begin{itemize}
    \item Charging and discharging power rating: 0.5 MW
    \item Energy capacity: 1 MWh
    \item Round-trip efficiency: 85\%
    \item Cycle life: 1000 cycles at 80\% cycle depth
    \item Shelf life: 5 years
    \item Pack cost: \$200/kWh
    \item Warranty coverage: until 80\% SoH
    \item True end-of-life: 75\%--50\% SoH
\end{itemize}

These cells have a near-quadratic cycle depth degradation stress function $\Phi$ that calculates the incremental capacity loss caused by a cycle of depth $u$~\cite{laresgoiti2015modeling}:
\begin{align}
    \Phi(u) = (3.14\text{E-4})u\up{2.03}\,.
    \label{Eq:DoD}
\end{align}
We do not consider power or C-rate as a degradation factor nor do we model the impact of degradation of the battery's rated power, because the assumed battery energy storage has a 2-hour duration corresponding to a maximum C-rate of 0.5, therefore impacts of C-rate magnitude over the battery degradation rate is negligible~\cite{saxena2019accelerated}. \rev{We also assume the battery has constant round-trip efficiency since at 0.5 C-rate lithium-ion batteries have stable voltage curve and impedance~\cite{chen2005thermal}}.

We use a rainflow-based degradation model to calculate the cycle degradation, which is represented as
\begin{align}
    D\ud{cyc} = \sum_{i=1}^{|\bm{u}|}\Phi(u_i),\;\bm{u} = \mathrm{Rainflow}(\bm{p}_n)
\end{align}
in which the total cycle degradation is calculated as the sum of cycle degradation caused by each cycle, and cycles are identified using the Rainflow algorithm using the storage operation profile. 

The daily calendar aging coefficient assuming a 5 year shelf life~\cite{bloom2010differential} is listed as
\begin{align}
    D\ud{cal} = 0.2/1825.
\end{align}

\revv{In all case studies, the battery degradation value function is discretized with a 1\% SoH granularity. For example, for a battery with 80\% EoL, the number of discretization segments $I$ as in \eqref{eq:pwl} is 21. }

\subsection{Battery Resale Value}

To better compare the value of a battery in grid services with its manufacturing cost, we assume a battery has a resale value while having a valid manufacturer warranty. For example, the battery can be sold to provide back-up power supply which is a much lower demanding usage scenario. We assume a battery has a resale value of \$200/kWh and has proratable resale values based on the remaining duration of the warranty and the remaining SoH. The resale value is described using the following formulation
\begin{align}
    S(E) = \begin{cases}
    \text{(\$200/kWh)}\Big(\frac{E-0.8E^0}{0.2E\up{0}}\Big)\Big(\frac{E}{E\up{0}}\Big) \quad & \text{if $E \geq 0.8E\up{0}$} \\
    0 & \text{else}
    \end{cases}
\end{align}
in which $S(E)$ models the resale value of the battery. \rev{We further assume the battery operator will optimize its battery operation considering the resale value, and update the battery degradation value as
\begin{align}
    V_n(E) = \max\{V_{n+1}(E), S(E)\}
\end{align}
meaning the battery operator should not use the battery to perform grid services at all as the income does not justify the reduction of the resale value. Instead, the battery operator should sell the battery cells at the resale price for a better return. }

\begin{figure}[t]%
	\centering
	\subfloat[New York State Price Zones.]{
		\includegraphics[trim = 0mm 0mm 0mm 0mm, clip, width = .85\columnwidth]{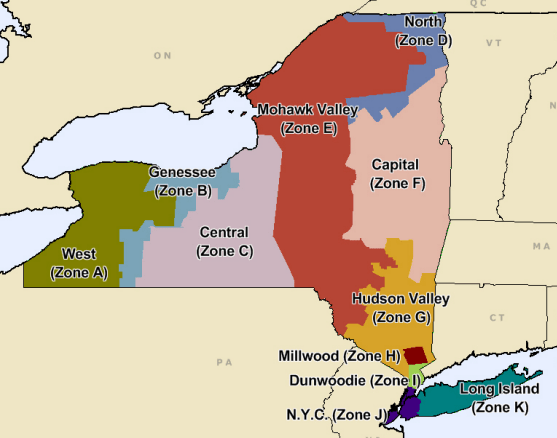}
		\label{fig:nyc1}%
	}
	\\
	\subfloat[30 day moving average price.]{
		\includegraphics[trim = 5mm 0mm 10mm 0mm, clip, width = .95\columnwidth]{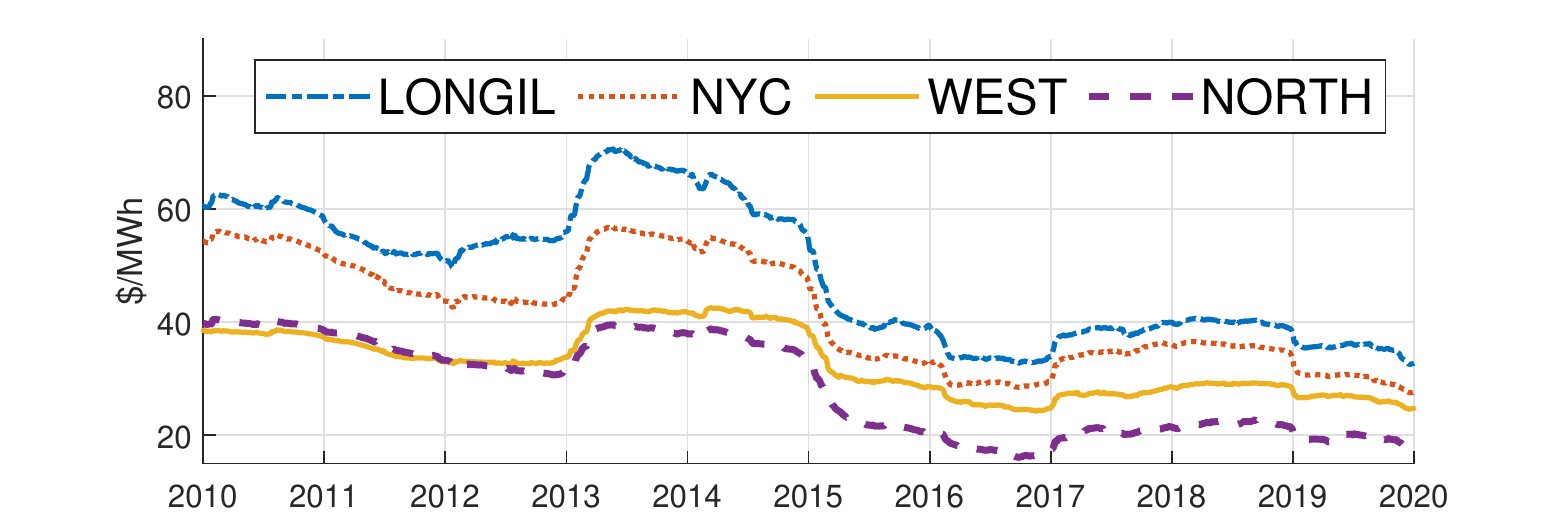}
		\label{fig:nyc2}%
	}
	\\
	\subfloat[30 day moving average daily price deviations.]{
		\includegraphics[trim = 5mm 0mm 10mm 0mm, clip, width = .95\columnwidth]{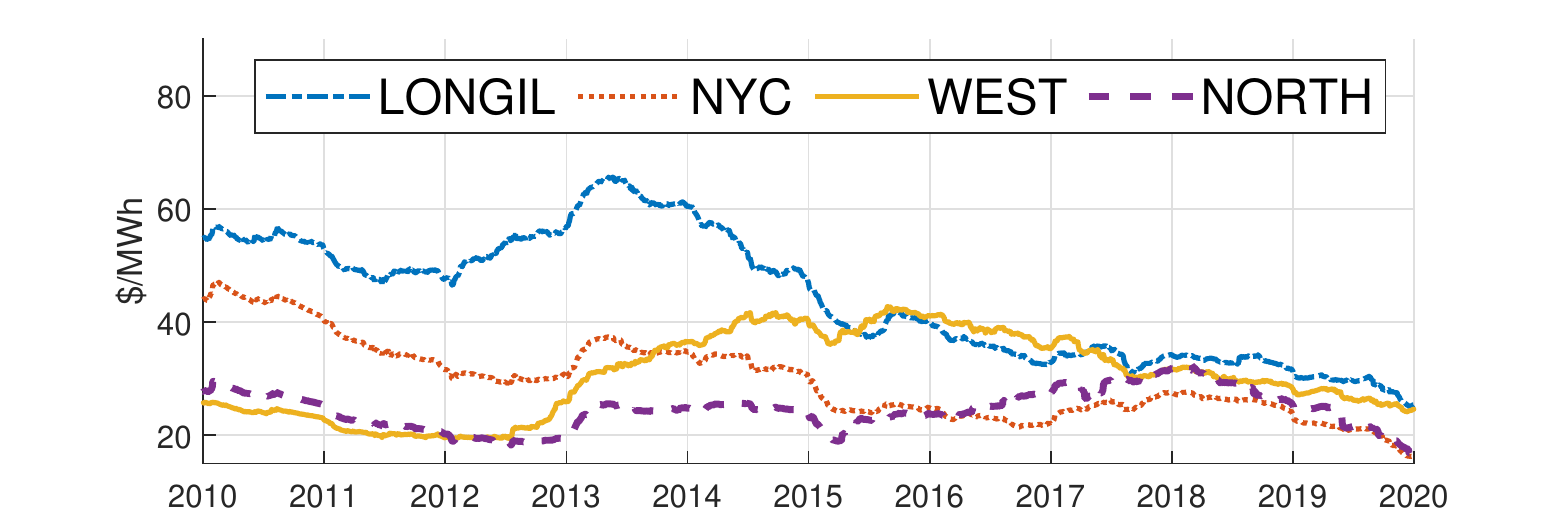}
		\label{fig:nyc3}%
	}
  \caption{Location of price zones and the historical price data of the selected four zones.}
    \label{fig:nyc}
\end{figure}

\subsection{Energy Price Arbitrage Case Study Settings}

We simulate a 10 year price arbitrage battery project from 2010 to 2020. 
This study assumes the battery energy storage can optimally arbitrage real-time price differences, which can either be accomplished by high-performance real-time price forecast model, or with properly designed market bidding models that will soon be available across all wholesale electricity markets in the U.S. according to FERC Order 841~\cite{sakti2018review}.

In price arbitrage applications the utilization function $O(\cdot)$ becomes
\begin{align}
    O_n(\bm{p}_n) = \sum_{t=1}^T M\lambda_np_{n,t}
\end{align}
where $M$ is the duration of each time step, which is 5-minute in the considered real-time arbitrage case study, and $\lambda_n$ is the market price. To avoid simultaneous charging and discharging during the occurrence of negative prices~\cite{xu2020lagrangian}, the discharge power constraint in \eqref{eq:dbv_c4} is slightly modified to enforce that the battery can only discharge when price is non-negative
\begin{align}
    0 \leq p_{n,t} \leq P\mathbf{1}_{\{\lambda_{n,t}\geq 0\}}
\end{align}
where $\mathbf{1}_{\{\lambda_{n,t}\geq 0\}}$ is the indicator function which returns 1 if $\lambda_{n,t}\geq 0$ and 0 otherwise. 
\rev{The arbitrage case study uses Algorithm 1 and incorporates the nonlinear cycle degradation model using a piecewise linear approximation approach with 10 segments, which achieves less than 1\% of approximation error~\cite{xu2017factoring}.} This piecewise linearization cycle degradation model with a total number of $J$ segments assuming a constant battery efficiency $\eta$ is summarized as
\begin{subequations}
\begin{align}
D_{cyc} &= \sum_{t=1}^T\sum_{j=1}^J \delta_jMp^{d}_{n,t,j} \\
p_{n,t} &= \sum_{j=1}^J p^{d}_{n,t,j} - p^{c}_{n,t,j}, \; p^{d}_{n,t,j}\geq 0,\; p^{c}_{n,t,j}\geq 0\\
e_{n,t,j} &= e_{n,t-1,j} + Mp^{c}_{n,t,j}\eta - Mp^{d}_{n,t,j}/\eta \\
0 &\leq  e_{n,t,j} \leq E_j
\end{align}
where $\delta_j$ is the incremental cycle degradation rate from the linearization segment $j$ and is calculated from the nonlinear cycle stress function as
\begin{align}
    \delta_j = \frac{J}{\eta E_n}\Big[\Phi(\frac{j}{J}) - \Phi(\frac{j-1}{J})\Big] 
\end{align}
and $p^{d}_{n,t,j}$ and $p^{c}_{n,t,j}$ and the discharge and charge power associated with each cycle depth linearization segment, $e_{n,t,j}$ is the auxiliary energy level and $E_j$ is the maximum energy level associated with each cycle depth segment.
\end{subequations}

We consider four price zones in NYISO: WEST (ZoneA), NORTH (ZoneD), NYC (Zone J), and LONGIL (Zone K). These four zones have distinct price behaviors due to congested transmission corridors and uneven resource distribution in New York State~\cite{patton20162014}. The location of the four zones, their average market prices and daily deviations, are demonstrated in Fig.~\ref{fig:nyc}. 
The model is implemented in Julia and solved using Gurobi on a personal computer. The solution time for finishing a 10-year valuation (3652 days each with 288 periods) took around one hour.

\subsection{\rev{Frequency Regulation  Case Study Settings}}


\begin{figure*}[t]%
    \vspace{-10mm}
	\centering
	\subfloat[Battery value in LONGIL including prorated resale values.]{
		\includegraphics[trim = 30mm 5mm 30mm 10mm, clip, width = .9\columnwidth]{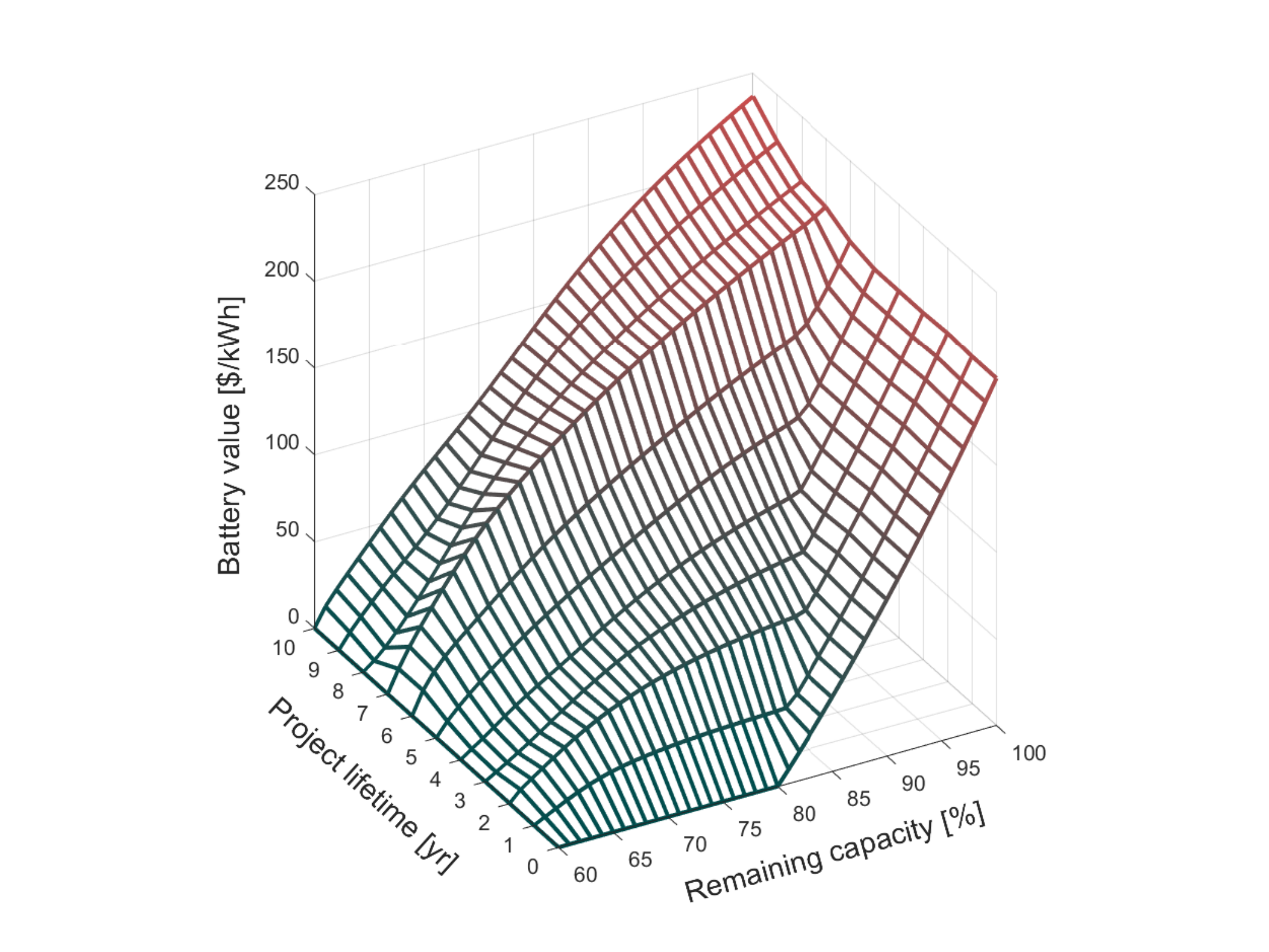}
		\label{fig:rev1}%
	}
	\subfloat[Surplus battery value in LONGIL.]{
		\includegraphics[trim = 30mm 5mm 30mm 10mm, clip, width = .9\columnwidth]{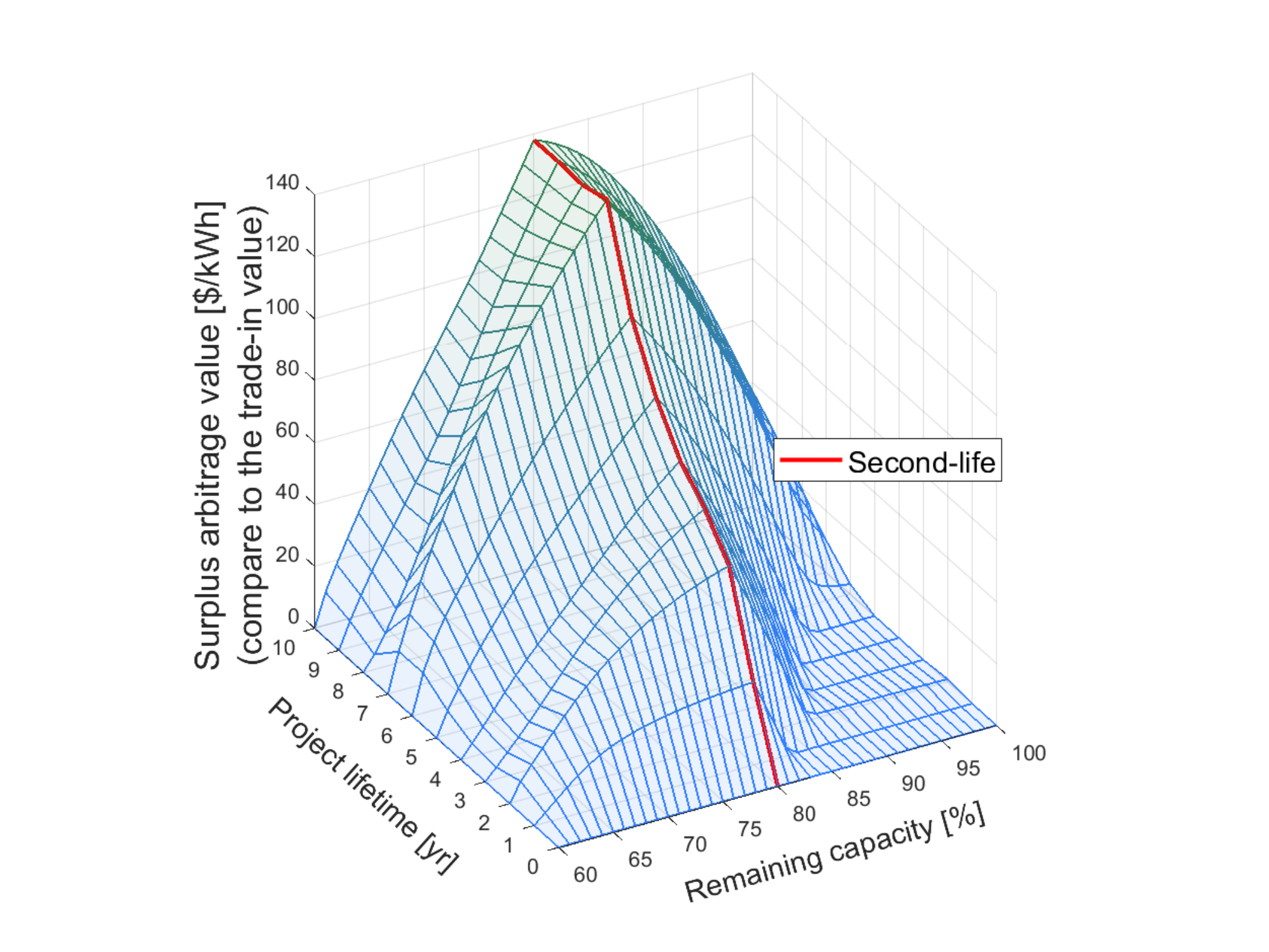}
		\label{fig:rev2}%
	}
	\\
	\subfloat[Incremental cost per EFC in LONGIL.]{
		\includegraphics[trim = 30mm 5mm 30mm 10mm, clip, width = .9\columnwidth]{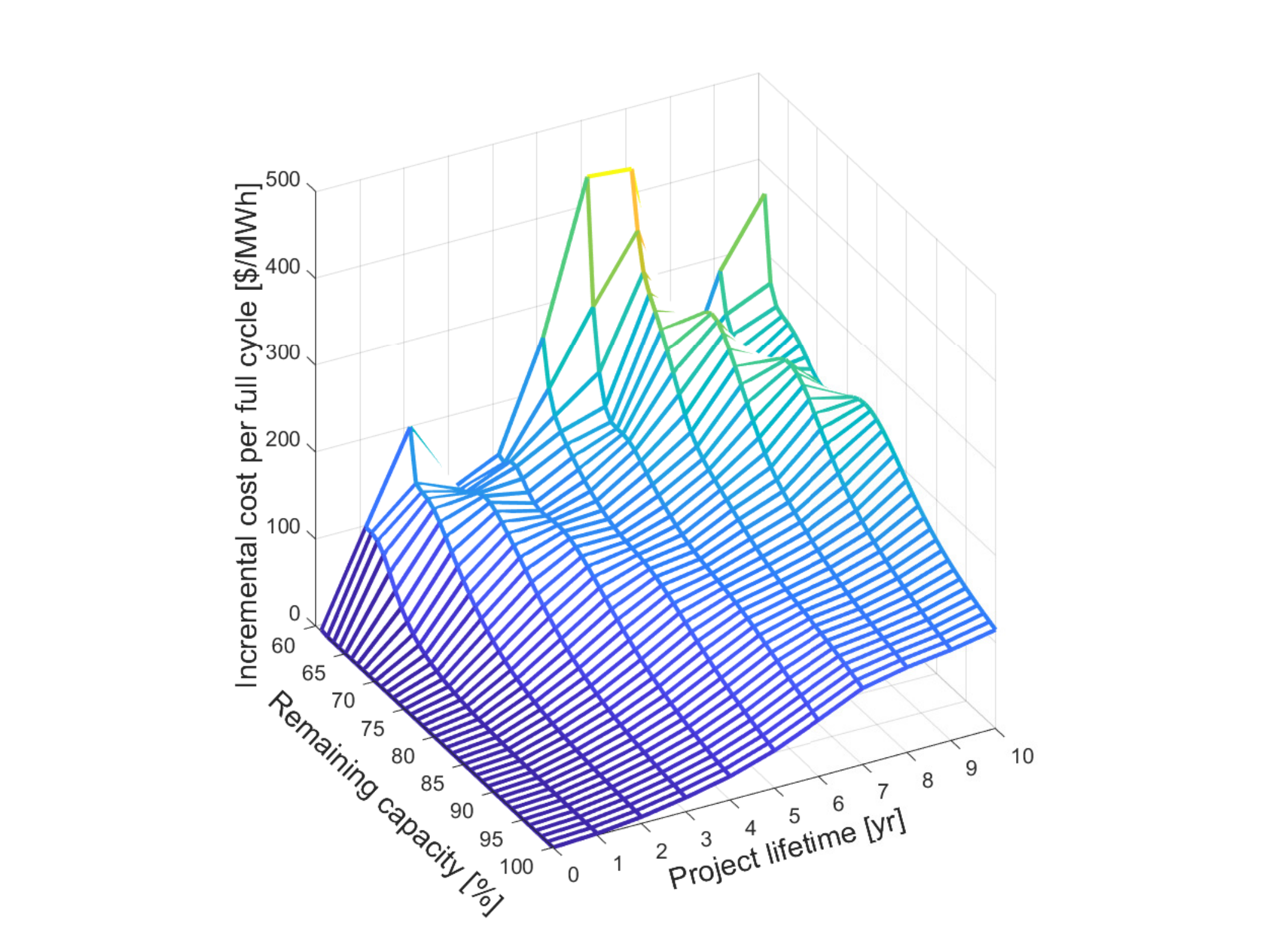}
		\label{fig:rev4}%
	}
	\subfloat[Comparison of surplus battery values.]{
		\includegraphics[trim = 0mm 0mm 0mm 0mm, clip, width = .9\columnwidth]{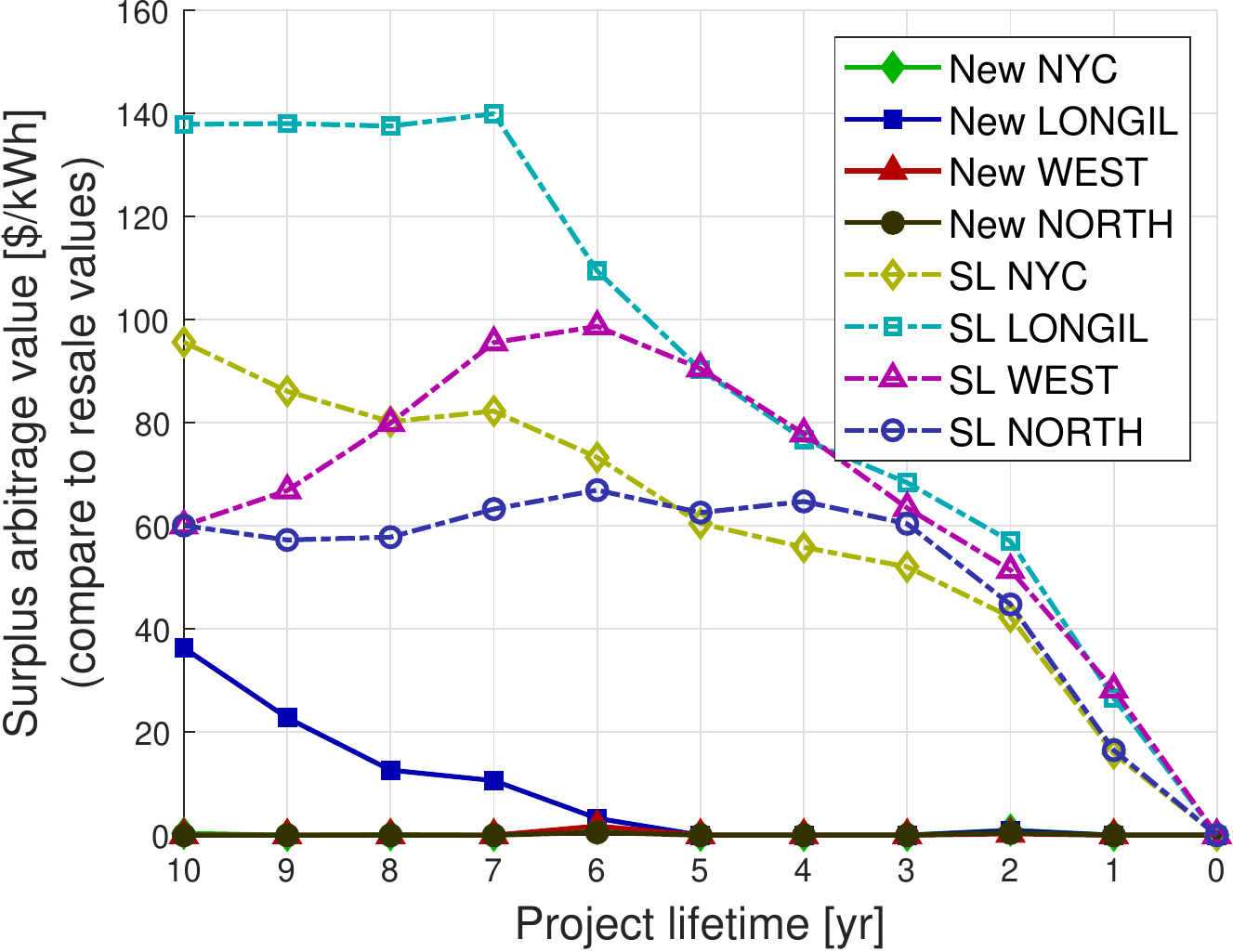}
		\label{fig:rev3}%
	}
  \caption{10 year arbitrage value analysis. (a) the total battery worth including resale value; (b) the surplus arbitrage value. (c) incremental cost per full cycle, resale value is not considered; (d) surplus value comparison of new and second-life (SL) batteries in the four NYISO price zones.}
    \label{fig:rev}
\end{figure*}

We use frequency regulation as the second case study to demonstrate the use of the dynamic valuation framework in a non-anticipatory control scenario (the battery responds to frequency regulation instruction signals without any knowledge about future signals). We consider the dynamic regulation service (RegD) in the PJM pay-for-performance frequency regulation market, one of the most popular markets for grid-scale battery projects. We use five years of frequency regulation price data from 2016 to 2020~\cite{pjm_reg}, to provide a quick context, the average daily payment for providing a 1~MW RegD regulation capacity is around \$550 per day.
On the operation side, we simulate the frequency regulation provision repetitively using the frequency regulation signal from 2020, due to PJM only provides historical signal data from 2020.

The battery's revenue in the PJM RegD regulation market over a provision period, which is an hour in PJM, can be represented as
\begin{align}
    \rho(\bm{p}_n, \bm{r}_n)\lambda_n B_n + \gamma V_{n+1}(E_{n+1})
\end{align}
where $\rho(\bm{p}_n, \bm{r}_n)$ calculates the performance score by comparing the regulation signal instruction $ \bm{r}_n$ and the battery response $\bm{p}_n$. This performance score is the average of an accuracy score (the relative absolute error), correlation score (the correlation coefficient), and a delay score (delay in responding to regulation instruction, battery usually achieves full delay score since the response delay is negligible compared to conventional generators), see the PJM market manual for details~\cite{pjm_manual}. $\lambda_n$ is the capacity clearing price and $B_n$ is the regulation capacity provided by the battery. $\gamma V_{n+1}(E_{n+1})$ is the opportunity value of the remaining battery lifetime as introduced in the dynamic valuation framework.

We assume the battery uses an optimal response controller, a result from our previous work~\cite{xu2018optimal}, as the control policy $G$ in our non-anticipatory dynamic valuation algorithm. This controller trades-off battery degradation cost and market performance penalties. This control policy limits the battery SoC swing within a range calculated as 
\begin{align}\label{eq:pol3}
    \hat{u} = \varphi^{-1}\Big(\frac{\eta^2+1}{\eta C_{i,n}}\pi\Big),\quad\pi = \frac{2}{3}\frac{\mu_{\lambda}}{\mu_{r}M} 
\end{align}
where $\varphi^{-1}(\cdot)$ is the inverse function of the derivative of the cycle stress function $\varphi(x) = \mathrm{d} \Phi(x)/\mathrm{d} x$. $\eta$ is the battery single-trip efficiency, $C_{i,n}$ is the marginal cost of degradation as shown in Algorithm 2, $\pi$ is calculated based on the expected regulation capacity price $\mu_{\lambda}$, the expected absolute sum of the regulation signal $\mu_{r}$, and the expected mileage ratio $M$.

The frequency regulation case study is implemented in Matlab and is performed according Algorithm 2, in which Step 2 performs the actual regulation simulation. 

\section{Results}

\subsection{Arbitrage Surplus Value Analysis}

We compare the battery value in conducting price arbitrage with a benchmark resale value in the four considered price zones in NYISO. Fig.~\ref{fig:rev} summarizes the result of this study. 
LONGIL achieves the highest valuation and is used as the demonstrative example. Fig.~\ref{fig:rev1} shows the battery value increases with the project duration and with the remaining capacity, while Fig.~\ref{fig:rev2} shows second-life batteries have the highest surplus valuation because as the battery ages, the battery resale value decreases at a faster rate than the arbitrage value, i.e., the marginal benefit of having a newer battery diminishes. This result can be better observed in Fig.~\ref{fig:rev4} that the incremental cost (decrements in the battery value) per EFC increases quickly as the battery ages, reaching up to a 10-time difference between new and near-EoL batteries. Value of batteries closer to EoL become very sensitive to price deviations as the remaining arbitrage window is very short, and gradually become smoother with higher SoH.

LONGIL is the only zone in which a new battery could be profitable due to its high price deviations, while Fig.~\ref{fig:rev3} shows that surplus values in the other three zones are zero most of the time. \rev{The only exceptions, in which new batteries are profitable in performing price arbitrage in price zones other than LONGIL, are observed between year 7 and year 6 due to the higher price average and deviations between year 2013 and 2014 as shown in Fig.~\ref{fig:nyc}.} However, second-life batteries achieved more than \$100/kWh valuation when the project has at least 5 years remaining, while in LONGIL second-life battery value could reach higher than \$200/kWh, surpassing the resale value of a new battery. These result concludes second-life batteries are competitive choices for conducting price arbitrage, and provides much higher surplus value to new batteries.

\subsection{\rev{Regulation Surplus Value Analysis}}

\begin{figure}[t]%
\centering
		\includegraphics[trim = 30mm 5mm 30mm 10mm, clip, width = .9\columnwidth]{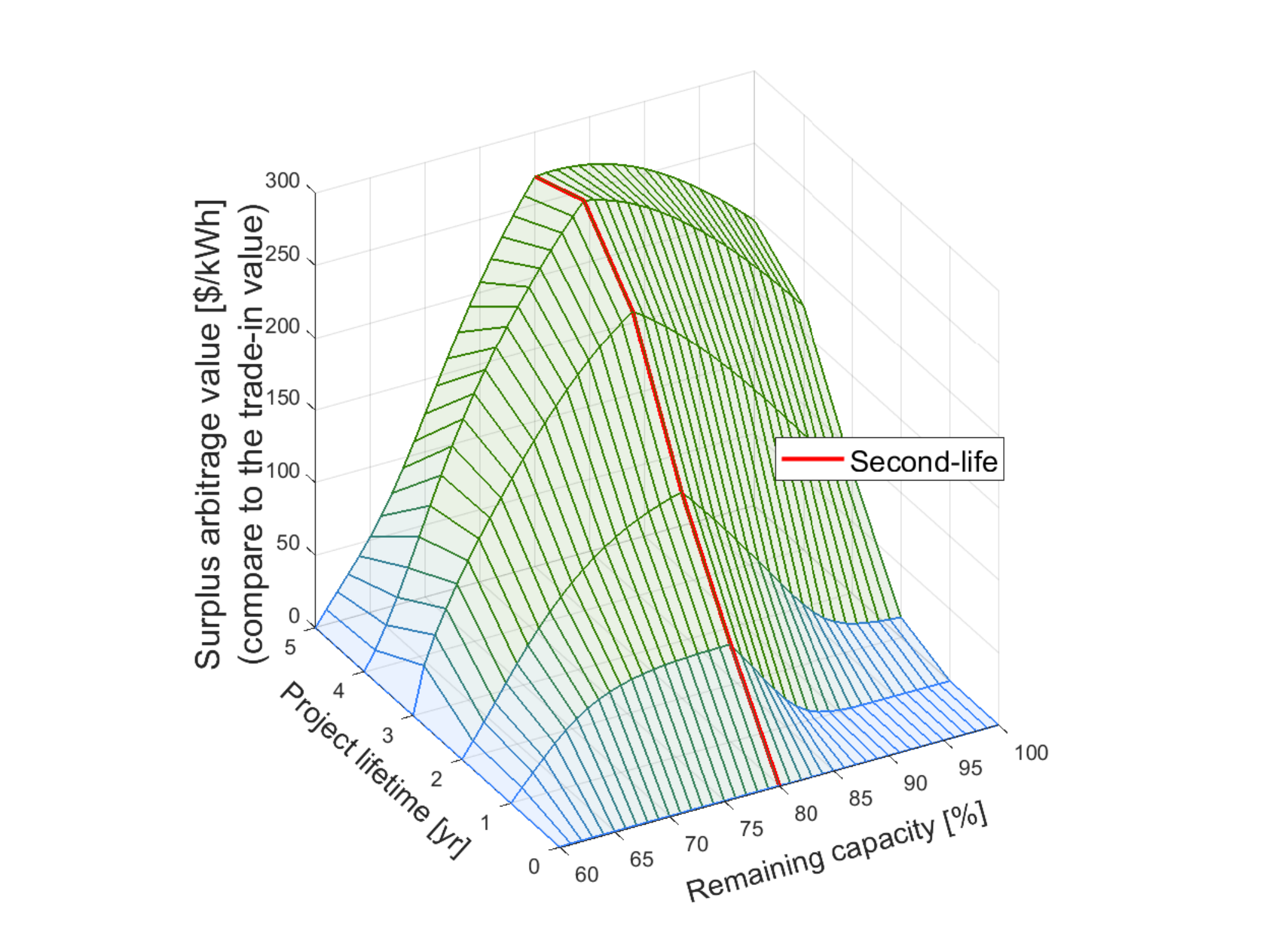}
  \caption{Surplus battery value in providing frequency regulation in PJM.}
    \label{fig:reg}
\end{figure}

Fig.~\ref{fig:reg} shows the surplus valuation result of the same battery in providing frequency regulation in PJM using a degradation-aware bidding and control strategy~\cite{xu2018optimal}. The battery earns significantly higher profit in providing frequency regulation and the surplus value is at least doubled compared to the arbitrage value from the LONGIL zone. Also note that the regulation case study only covers five years due to limited data availability, so at the five year project period the regulation case makes about more than three time more revenue compared to the arbitrage case. 

The regulation case study results also show that 80\% SoH provides the best surplus value. Yet, the differences among 80\% to 100\% is not significant compared to the arbitrage case. The surplus value of an 80\% SoH battery is about 20\% higher than that of a new battery. Also, due to the limited price data availability, a five-year project duration does not seem to provide full utilization of a new battery, as the figure shows the battery value kept rising from year 0 to year 5. Comparably, the value of an 80\% SoH battery saturated when the project has more than 4 years left. These results show that new batteries may be better candidates for building frequency regulation projects, especially considering the additional cost of battery recycling and purchase delivery. 

\subsection{The Value of Second-life Batteries}\label{sec:sl}

\begin{figure*}[t]%
	\centering
	\subfloat[NYISO NYC arbitrage]{
		\includegraphics[trim = 0mm 00mm 0mm 0mm, clip, width = .55\columnwidth]{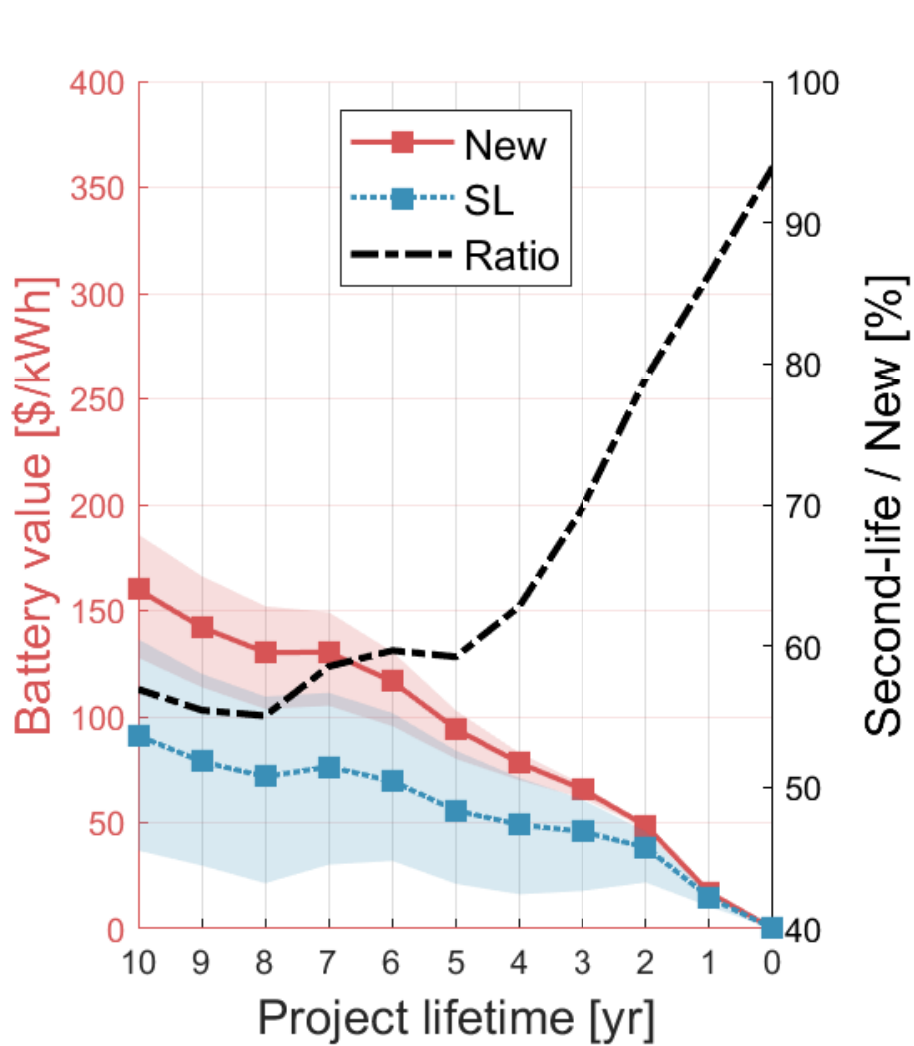}
		\label{fig:sl1}%
	}
	\subfloat[NYISO LONGIL arbitrage]{
		\includegraphics[trim = 0mm 00mm 0mm 0mm, clip, width = .55\columnwidth]{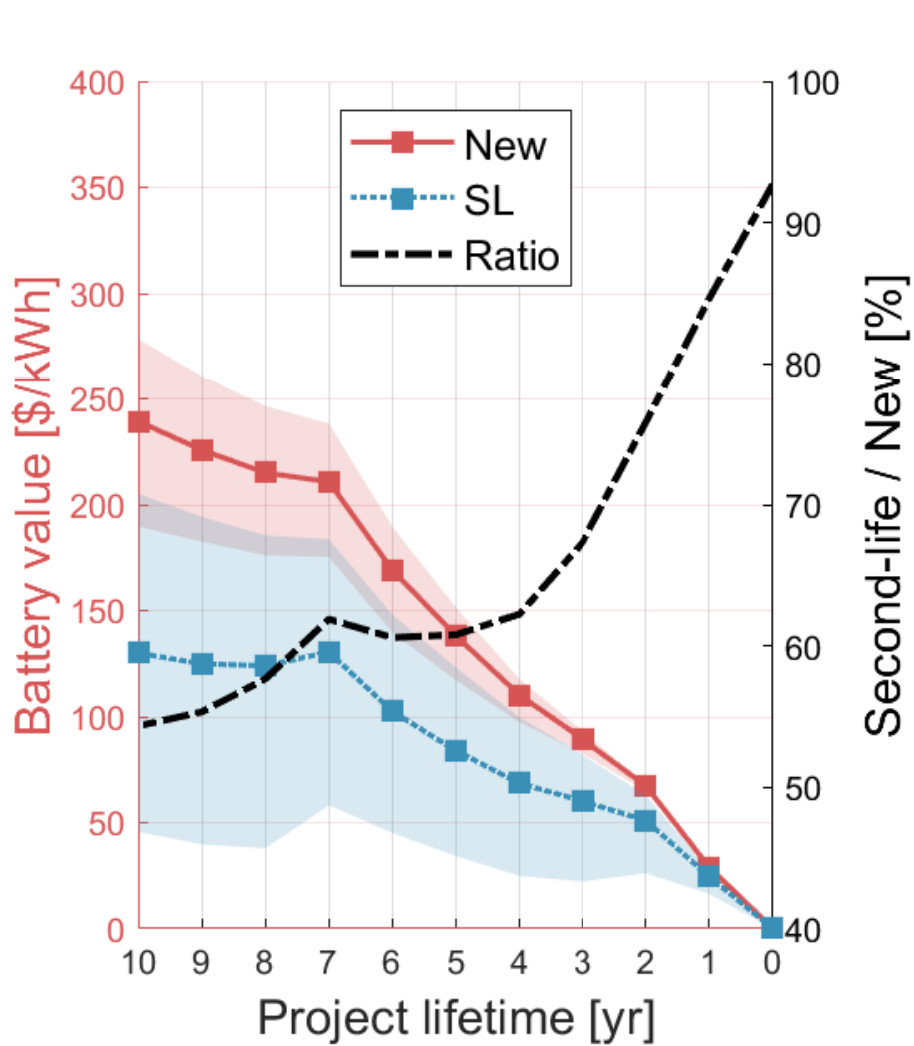}
		\label{fig:sl2}%
	}
	\subfloat[PJM RegD w/o degradation control]{
		\includegraphics[trim = 0mm 00mm 0mm 0mm, clip, width = .55\columnwidth]{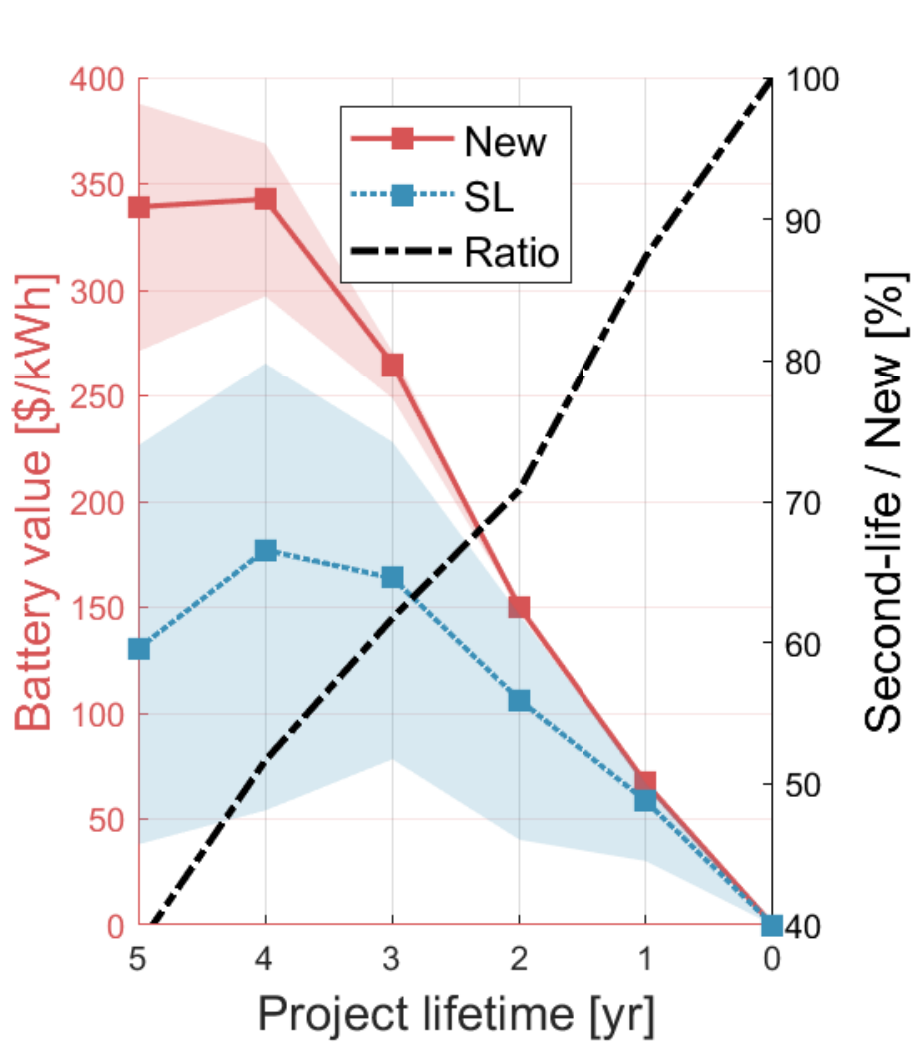}
		\label{fig:sl5}%
	}	\\
	\subfloat[NYISO NORTH arbitrage]{
		\includegraphics[trim = 0mm 0mm 0mm 0mm, clip, width = .55\columnwidth]{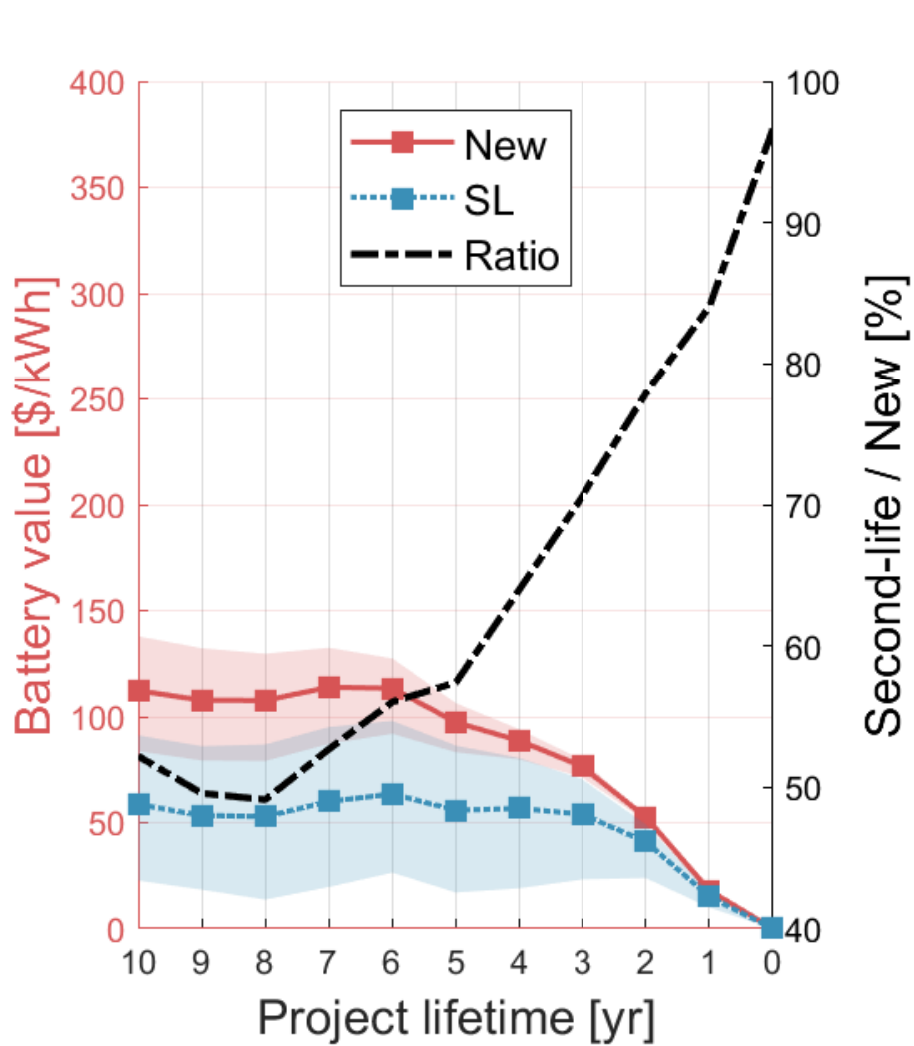}
		\label{fig:sl3}%
	}
	\subfloat[NYISO WEST arbitrage]{
		\includegraphics[trim = 0mm 0mm 0mm 0mm, clip, width = .55\columnwidth]{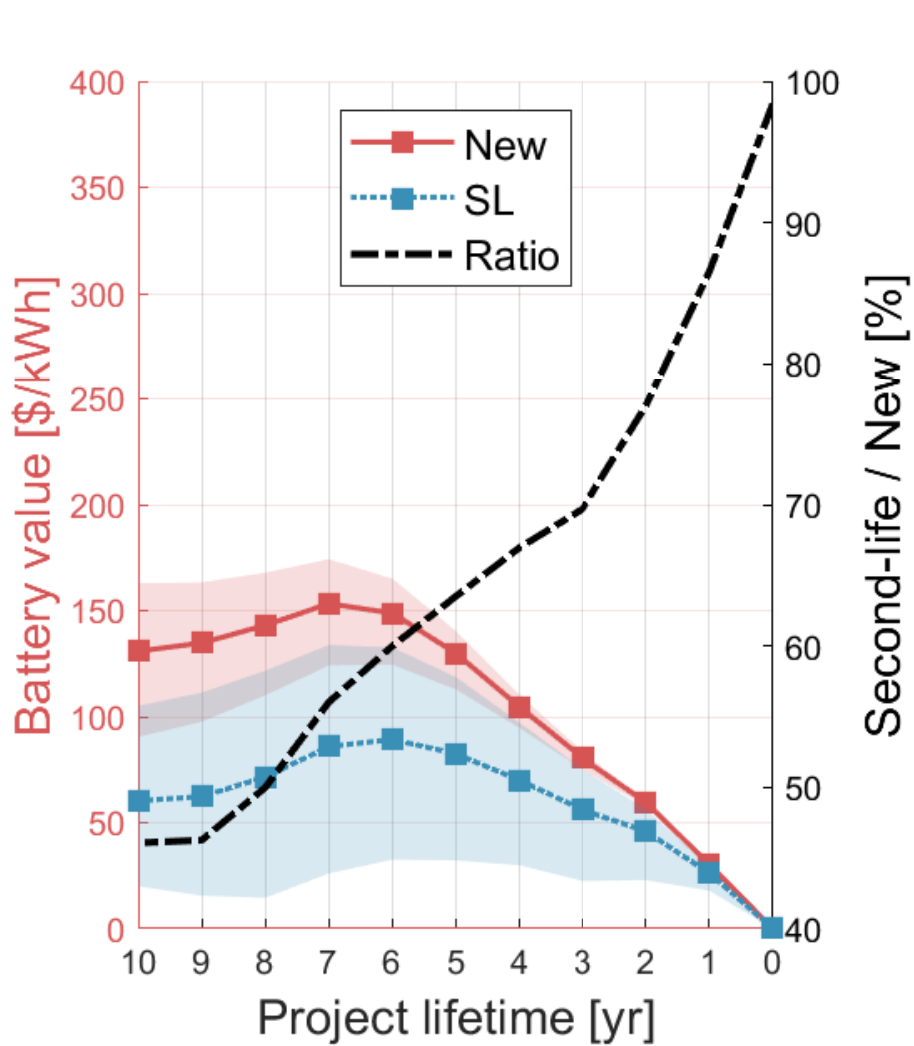}
		\label{fig:sl4}%
	}
	\subfloat[PJM RegD with degradation control]{
		\includegraphics[trim = 0mm 00mm 0mm 0mm, clip, width = .55\columnwidth]{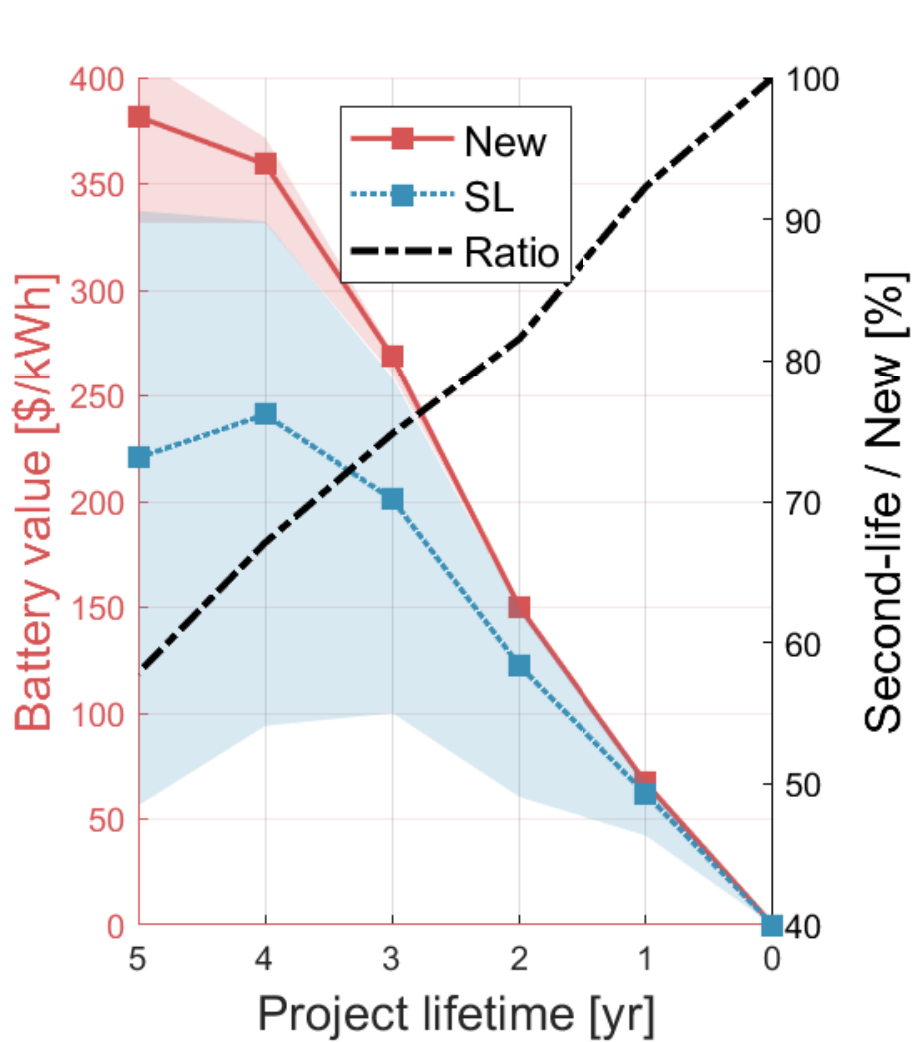}
		\label{fig:sl6}%
	}
  \caption{Comparison of the battery lifetime value between new and second-life (80\% state of life) batteries in arbitrage and frequency regulation, shaded areas represent the value range of EoL between 75\% to 50\% SoH.}
    \label{fig:sl}
\end{figure*}

In practice there is no guarantee over the battery end-of-life when the warranty expired, thus the battery could become unusable at any time and must be terminated from service. In this section, we assume the battery has an uniformly distributed end-of-life from 75\% to 50\% SoH, and compare the expected value of new and second-life batteries given this EoL distribution. To implement this uniform distribution assumption of EoL, we assume the battery end-of-life has 6 possibilities: 50\%, 55\%, 60\%, 65\%, 70\%, 75\%, each with one-sixth of possibility. And perform a battery valuation over each of the EoL scenario, the valuation results are demonstrated in Fig.~\ref{fig:sl}. At each valuation point, we do a weighted sum of the battery value from the six EoL scenarios to calculate the average new and second-life battery value, and the ratio between them.

Fig.~\ref{fig:sl} shows the average value of new and second-life batteries in performing price arbitrage in four considered price zones in NYISO, and two results of providing frequency regulation in PJM: one not using any degradation-aware control so the 0.5MW/1MWh battery always provides 0.5MW regulation capacity; the other case adopts a degradation-aware bidding and control strategy. Note that we do not consider the resale value in this study to provide a fair comparison between new and second-life batteries. The valuation result for new batteries covers both the warranty period (above 80\% SoH) and the second-life period (80\% SoH to EoL). In general, over a 10-year arbitrage project, second-life battery values are much more sensitive to EoL, while new battery values are more stable. 

In the frequency regulation case, degradation-aware strategies improve the battery value, but not a lot in new batteries. The primary reason is that the battery considered in this case study has a two-hour duration, thus providing frequency regulation will not lead to many deep cycles that significantly reduces battery lifetime. The use of degradation-aware control is more obvious in second-life batteries, improving about 50\% of the lifetime value.

A key takeaway is that  ratios between second-life and new batteries are very similar in the four arbitrage cases and the regulation case with degradation-aware control, despite these cases have distinct price behaviors and battery values. As the project duration shortens, the worth of second-life batteries increases steadily from around 60\% to around 95\% of new batteries, showing very promising economic value especially considering the average lifetime energy throughput of second-life batteries is only around 42\% of new batteries.

\section{Conclusion}

In this study, we developed a dynamic valuation framework that finds the financial worth of batteries participating in electricity markets given the remaining project duration and battery lifetime, and demonstrates results from price arbitrage in NYISO zones and frequency regulation in PJM using  historical data. \rev{In addition to these two application scenarios, the framework provides wide applicability. It is applicable to perform long-term battery capacity and lifetime valuation in many different battery energy storage services, including stochastic price arbitrage, peak shaving, and stacked services.}

This paper also performed a value analysis for second-life batteries, considering the uncertainties of their end-of-life due to the absence of manufacture warranties. Results show that values of second-life batteries follow a stable ratio to new batteries in various distinct market scenarios. This ratio steadily increases from around 60\% at the start of a 10-year project to around 95\% at the project end, making second-life batteries a much economical choice for conducting price arbitrage in electricity markets.

With the battery economy continuing to boom,  manufacturing costs of batteries will no longer be a limiting factor for the grid-scale applications, and second-life batteries retired from electric vehicles pose vast opportunities to be re-purposed for grid applications. Results from this study provide vital reference over the worth of new and used batteries in electricity markets, concluding that second-life batteries are superior choices for grid applications.  However, this study neglected costs associated with utilizing second-life batteries, including logistics and assembling. Nevertheless, we hope these valuation results will provide incentive industry stakeholders to develop new standards to facilitate the reuse of second-life batteries, and improve the public awareness over the value of used batteries.

\bibliographystyle{IEEEtran}	
\bibliography{IEEEabrv,sample}		

\begin{thebibliography}{10}
\providecommand{\url}[1]{#1}
\csname url@samestyle\endcsname
\providecommand{\newblock}{\relax}
\providecommand{\bibinfo}[2]{#2}
\providecommand{\BIBentrySTDinterwordspacing}{\spaceskip=0pt\relax}
\providecommand{\BIBentryALTinterwordstretchfactor}{4}
\providecommand{\BIBentryALTinterwordspacing}{\spaceskip=\fontdimen2\font plus
\BIBentryALTinterwordstretchfactor\fontdimen3\font minus
  \fontdimen4\font\relax}
\providecommand{\BIBforeignlanguage}[2]{{%
\expandafter\ifx\csname l@#1\endcsname\relax
\typeout{** WARNING: IEEEtran.bst: No hyphenation pattern has been}%
\typeout{** loaded for the language `#1'. Using the pattern for}%
\typeout{** the default language instead.}%
\else
\language=\csname l@#1\endcsname
\fi
#2}}
\providecommand{\BIBdecl}{\relax}
\BIBdecl

\bibitem{blochbreakthrough}
C.~Bloch, J.~Newcomb, S.~Shiledar, and M.~Tyson, ``Breakthrough batteries,''
  \emph{Powering the Era of Clean Electrification}, 2020.

\bibitem{david2001market}
A.~K. David and F.~Wen, ``Market power in electricity supply,'' \emph{IEEE
  Transactions on energy conversion}, vol.~16, no.~4, pp. 352--360, 2001.

\bibitem{vetter2005ageing}
J.~Vetter, P.~Nov{\'a}k, M.~R. Wagner, C.~Veit, K.-C. M{\"o}ller, J.~Besenhard,
  M.~Winter, M.~Wohlfahrt-Mehrens, C.~Vogler, and A.~Hammouche, ``Ageing
  mechanisms in lithium-ion batteries,'' \emph{Journal of power sources}, vol.
  147, no. 1-2, pp. 269--281, 2005.

\bibitem{preger2020degradation}
Y.~Preger, H.~M. Barkholtz, A.~Fresquez, D.~L. Campbell, B.~W. Juba,
  J.~Rom{\`a}n-Kustas, S.~R. Ferreira, and B.~Chalamala, ``Degradation of
  commercial lithium-ion cells as a function of chemistry and cycling
  conditions,'' \emph{Journal of The Electrochemical Society}, vol. 167,
  no.~12, p. 120532, 2020.

\bibitem{kumar2012power}
N.~Kumar, P.~Besuner, S.~Lefton, D.~Agan, and D.~Hilleman, ``Power plant
  cycling costs,'' National Renewable Energy Lab.(NREL), Golden, CO (United
  States), Tech. Rep., 2012.

\bibitem{xu2016comparison}
B.~Xu, Y.~Dvorkin, D.~S. Kirschen, C.~A. Silva-Monroy, and J.-P. Watson, ``A
  comparison of policies on the participation of storage in us frequency
  regulation markets,'' in \emph{2016 IEEE Power and Energy Society General
  Meeting (PESGM)}.\hskip 1em plus 0.5em minus 0.4em\relax IEEE, 2016, pp.
  1--5.

\bibitem{zakeri2015electrical}
B.~Zakeri and S.~Syri, ``Electrical energy storage systems: A comparative life
  cycle cost analysis,'' \emph{Renewable and sustainable energy reviews},
  vol.~42, pp. 569--596, 2015.

\bibitem{doe_ev}
{Department of Energy}, ``Electric car safety, maintenance, and battery life,''
  2020, [Available Online]
  \url{https://www.energy.gov/eere/electricvehicles/electric-car-safety-maintenance-and-battery-life}.

\bibitem{xu2017scalable}
B.~Xu, Y.~Wang, Y.~Dvorkin, R.~Fern{\'a}ndez-Blanco, C.~A. Silva-Monroy, J.-P.
  Watson, and D.~S. Kirschen, ``Scalable planning for energy storage in energy
  and reserve markets,'' \emph{IEEE Transactions on Power systems}, vol.~32,
  no.~6, pp. 4515--4527, 2017.

\bibitem{wang2019power}
Y.~Wang, S.~J. Moura, S.~G. Advani, and A.~K. Prasad, ``Power management system
  for a fuel cell/battery hybrid vehicle incorporating fuel cell and battery
  degradation,'' \emph{International Journal of Hydrogen Energy}, vol.~44,
  no.~16, pp. 8479--8492, 2019.

\bibitem{xu2017factoring}
B.~Xu, J.~Zhao, T.~Zheng, E.~Litvinov, and D.~S. Kirschen, ``Factoring the
  cycle aging cost of batteries participating in electricity markets,''
  \emph{IEEE Transactions on Power Systems}, vol.~33, no.~2, pp. 2248--2259,
  2017.

\bibitem{shi2017using}
Y.~Shi, B.~Xu, D.~Wang, and B.~Zhang, ``Using battery storage for peak shaving
  and frequency regulation: Joint optimization for superlinear gains,''
  \emph{IEEE Transactions on Power Systems}, vol.~33, no.~3, pp. 2882--2894,
  2017.

\bibitem{he2015optimal}
G.~He, Q.~Chen, C.~Kang, P.~Pinson, and Q.~Xia, ``Optimal bidding strategy of
  battery storage in power markets considering performance-based regulation and
  battery cycle life,'' \emph{IEEE Transactions on Smart Grid}, vol.~7, no.~5,
  pp. 2359--2367, 2015.

\bibitem{reniers2018improving}
J.~M. Reniers, G.~Mulder, S.~Ober-Bl{\"o}baum, and D.~A. Howey, ``Improving
  optimal control of grid-connected lithium-ion batteries through more accurate
  battery and degradation modelling,'' \emph{Journal of Power Sources}, vol.
  379, pp. 91--102, 2018.

\bibitem{mohsenian2015optimal}
H.~Mohsenian-Rad, ``Optimal bidding, scheduling, and deployment of battery
  systems in california day-ahead energy market,'' \emph{IEEE Transactions on
  Power Systems}, vol.~31, no.~1, pp. 442--453, 2015.

\bibitem{he2018intertemporal}
G.~He, Q.~Chen, P.~Moutis, S.~Kar, and J.~F. Whitacre, ``An intertemporal
  decision framework for electrochemical energy storage management,''
  \emph{Nature Energy}, vol.~3, no.~5, pp. 404--412, 2018.

\bibitem{he2020power}
G.~He, S.~Kar, J.~Mohammadi, P.~Moutis, and J.~Whitacre, ``Power system
  dispatch with marginal degradation cost of battery storage,'' \emph{IEEE
  Transactions on Power Systems}, 2020.

\bibitem{birkl2017degradation}
C.~R. Birkl, E.~McTurk, S.~Zekoll, F.~H. Richter, M.~R. Roberts, P.~G. Bruce,
  and D.~A. Howey, ``Degradation diagnostics for commercial lithium-ion cells
  tested at- 10° c,'' \emph{Journal of the Electrochemical Society}, vol. 164,
  no.~12, p. A2644, 2017.

\bibitem{bellman2015applied}
R.~E. Bellman and S.~E. Dreyfus, \emph{Applied dynamic programming}.\hskip 1em
  plus 0.5em minus 0.4em\relax Princeton university press, 2015.

\bibitem{kirschen2018fundamentals}
D.~S. Kirschen and G.~Strbac, \emph{Fundamentals of power system
  economics}.\hskip 1em plus 0.5em minus 0.4em\relax John Wiley \& Sons, 2018.

\bibitem{krishnamurthy2017energy}
D.~Krishnamurthy, C.~Uckun, Z.~Zhou, P.~R. Thimmapuram, and A.~Botterud,
  ``Energy storage arbitrage under day-ahead and real-time price uncertainty,''
  \emph{IEEE Transactions on Power Systems}, vol.~33, no.~1, pp. 84--93, 2017.

\bibitem{xu2018optimal}
B.~Xu, Y.~Shi, D.~S. Kirschen, and B.~Zhang, ``Optimal battery participation in
  frequency regulation markets,'' \emph{IEEE Transactions on Power Systems},
  vol.~33, no.~6, pp. 6715--6725, 2018.

\bibitem{bishop2013evaluating}
J.~D. Bishop, C.~J. Axon, D.~Bonilla, M.~Tran, D.~Banister, and M.~D.
  McCulloch, ``Evaluating the impact of v2g services on the degradation of
  batteries in phev and ev,'' \emph{Applied energy}, vol. 111, pp. 206--218,
  2013.

\bibitem{shi2018convex}
Y.~Shi, B.~Xu, Y.~Tan, and B.~Zhang, ``A convex cycle-based degradation model
  for battery energy storage planning and operation,'' in \emph{2018 Annual
  American Control Conference (ACC)}.\hskip 1em plus 0.5em minus 0.4em\relax
  IEEE, 2018, pp. 4590--4596.

\bibitem{ecker2014calendar}
M.~Ecker, N.~Nieto, S.~K{\"a}bitz, J.~Schmalstieg, H.~Blanke, A.~Warnecke, and
  D.~U. Sauer, ``Calendar and cycle life study of li (nimnco) o2-based 18650
  lithium-ion batteries,'' \emph{Journal of Power Sources}, vol. 248, pp.
  839--851, 2014.

\bibitem{laresgoiti2015modeling}
I.~Laresgoiti, S.~K{\"a}bitz, M.~Ecker, and D.~U. Sauer, ``Modeling mechanical
  degradation in lithium ion batteries during cycling: Solid electrolyte
  interphase fracture,'' \emph{Journal of Power Sources}, vol. 300, pp.
  112--122, 2015.

\bibitem{saxena2019accelerated}
S.~Saxena, Y.~Xing, D.~Kwon, and M.~Pecht, ``Accelerated degradation model for
  c-rate loading of lithium-ion batteries,'' \emph{International journal of
  electrical power \& energy systems}, vol. 107, pp. 438--445, 2019.

\bibitem{chen2005thermal}
S.~Chen, C.~Wan, and Y.~Wang, ``Thermal analysis of lithium-ion batteries,''
  \emph{Journal of power sources}, vol. 140, no.~1, pp. 111--124, 2005.

\bibitem{sakti2018review}
A.~Sakti, A.~Botterud, and F.~O’Sullivan, ``Review of wholesale markets and
  regulations for advanced energy storage services in the united states:
  Current status and path forward,'' \emph{Energy policy}, vol. 120, pp.
  569--579, 2018.

\bibitem{xu2020lagrangian}
B.~Xu, M.~Korp{\aa}s, A.~Botterud, and F.~O’Sullivan, ``A lagrangian policy
  for optimal energy storage control,'' in \emph{2020 American Control
  Conference (ACC)}.\hskip 1em plus 0.5em minus 0.4em\relax IEEE, 2020, pp.
  224--230.

\bibitem{patton20162014}
D.~B. Patton, P.~LeeVanSchaick, J.~Chen, and M.~M. Unit, ``2014 state of the
  market report for the new york iso markets,'' \emph{Potomac Economics}, 2016.

\bibitem{pjm_reg}
``Pjm ancillary services,'' [Available Online]
  \url{https://www.pjm.com/markets-and-operations/ancillary-services.aspx}.

\bibitem{pjm_manual}
``Pjm manual 12: Balancing operations,'' [Available Online]
  \url{https://www.pjm.com/~/media/documents/manuals/m12.ashx}.

\end{thebibliography}

\begin{IEEEbiographynophoto}{Bolun Xu}
(S'14-M'18) received B.S. degrees from Shanghai Jiaotong
University, Shanghai, China in 2011;  M.Sc degree from Swiss Federal Institute of Technology, Zurich, Switzerland in 2014; and Ph.D. degree from University of Washington, Seattle, U.S. in 2018; all from Electrical Engineering.

He is currently an assistant professor in Columbia University, Department of Earth and Environmental Engineering, with affiliation in Department of Electrical Engineering. His research interests include electricity markets, energy storage, power system optimization, and power system economics. 
\end{IEEEbiographynophoto}

\end{document}